\documentclass[10pt]{IEEEtran}
\usepackage{amsmath,amsfonts,amsthm,bm,cancel}
\usepackage[margin=0.5in]{geometry}
\usepackage{multicol}
\usepackage{caption}
\usepackage{graphicx}

\title{A Chance Constraint Predictive Control and Estimation Framework for Spacecraft Descent with Field Of View Constraints}

\author{Steven van Leeuwen \thanks{Steven van Leeuwen is with the Department of Aerospace Engineering, University of Michigan- Ann Arbor, MI.}}

\begin{document}
\maketitle
\begin{abstract}
Recent studies of optimization methods and GNC of spacecraft near small bodies focusing on descent, landing, rendezvous, etc., with key safety constraints such as line-of-sight conic zones and soft landings have shown promising results; this paper considers descent missions to an asteroid surface with a constraint that consists of an onboard camera and asteroid surface markers while using a stochastic convex MPC law. An undermodeled asteroid gravity and spacecraft technology inspired measurement model is established to develop the constraint. Then a computationally light stochastic Linear Quadratic MPC strategy is presented to keep the spacecraft in satisfactory field of view of the surface markers while trajectory tracking, employing chance based constraints and up-to-date estimation uncertainty from navigation. The estimation uncertainty giving rise to the tightened constraints is particularly addressed. Results suggest robust tracking performance across a variety of trajectories. \end{abstract}
\section{Introduction}
Control of spacecraft in space missions such as rendezvous, docking at a station, landing on an orbital body, etc. are increasingly gathering interest as space capabilities improve and space missions continue, such as the ongoing Hayabusa2 mission [6] which is an asteroid sample and return to earth mission. The treatment of constraints in these missions is critical, and MPC strategies have had their fair share of success in various simulated missions. Constraints employed include glide-slope constraints, line-of-sight conic zones or parabolic zone constraints, etc., as a function of state [1,3,8]. However with focus on these aspects in the control formulation, the measurement model that arrives at the state estimation and the accompanying noise is often overlooked or simplified. In missions the onboard camera can play an important role in navigation [13,14] and thus is the motivation for an optical field of view constraint that aims to keep certain asteroid surface markers in the field of view, as is the case in the Hayabusa2 mission.

In the particular setting of a spacecraft descending toward an asteroid surface, uncertainty in the asteroid gravity as well as in the measurement equipment adds to the compexity of the control problem and mission. MPC with chance constraints, a basic form of stochastic MPC, is an approach to deal with constraints and trajectory tracking when faced with uncertainty and is explored in this paper. The asteroid used here is 433 Eros and its gravity is calculated as a high-fidelity model based on spherical harmonic expansion for use in the plant, and a low-fidelity model based on a constant density ellipsoid and truncated high-fidelity model for use in the control, details regarding both gravity models are in [8]. Navigation is achieved through use of an Extended Kalman Filter (EKF), and the up-to-date estimation error is used at each time-step in the MPC formulation, as was done in [1]. 

This paper contributes a technology inspired measurement model that gives rise to the estimation uncertainty and field of view constraint, and an accompanying stochastic tube MPC control.

The paper is layed out as follows. Notations used are presented first. Then follows plant, control, measurement, and estimation model descriptions, with particular explanation given to the measurement model that allows for a field of view constraint. Section IV describes how the chance based constraint is built and the MPC formulation. Section V then gives the results of two simulated trajectories.
\section{Definitions}
The controlled descent of the spacecraft around the asteroid is referred to as the spacecraft mission. Define the inertial frame, centered at the origin of the asteroid, by subscript $i$. The asteroid is considered in constant rotation, described by the asteroid fixed frame with subscript $a$. This frame rotates with respect to the inertial frame. The spacecraft body frame is denoted by subscript $b$. As an example of two subscripts written together, $ba$ gives a quantity of interest in the asteroid fixed frame expressed in the spacecraft body frame, or described a transformation from $a$ to $b$ when attached to a matrix. Bolded notation denotes vectors, and is reserved for model parameters.

Define the following:
$\mathbf{a}_a = [\ddot{x}\ \ddot{y}\ \ddot{z}]^T$,
$\mathbf{v}_a = [\dot{x}\ \dot{y}\ \dot{z}]^T$,
$\mathbf{r}_a = [x\ y\ z]^T$,
$\mathbf{u} := \mathbf{u}_i = [u_x\ u_y\ u_z]^T$, with $u_x, u_y, u_z$ the control forces of the spacecraft per unit mass expressed in the inertial, and $\mathbf{F} = [F_{g,x}\ F_{g,y}\ F_{g,z}]^T$ which is the high fidelity model of gravity for the asteroid [3,8] in the inertial frame. $\hat{\mathbf{F}}$ is the low fidelity model of gravity. 

There is also $\mathbf{\Omega} = [0\ 0\ n_i]^T$, the rotation vector of the asteroid expressed in the inertial frame, $\mu$ the gravitational constant of the asteroid, and $\mathbf{\Theta} := \mathbf{\Theta}_{bi} = [\phi\ \theta\ \psi]^T$, where $\phi,\theta,\psi$ are 321 Euler angles, with $\psi$ representing the first transformation from the inertial frame to the spacecraft Euler 1 frame. $\boldsymbol{\omega}_b := \boldsymbol{\omega}_{bi} = [\omega_x\ \omega_y\ \omega_z]^T$ are the angular velocities of the spacecraft. $\mathbf{J}$ is the inertia matrix of the spacecraft. The control moments of the spacecraft, not per unit mass, are
$\mathbf{M} := \mathbf{M}_{bi} = [M_x\ M_y\ M_z]^T$. 

In section IV the subscript $a|b$ denotes a quantity computed $a$ time steps into the future starting from time $b$. $\oplus,\ominus$ denote Minkowski addition and Minkowski subtraction (P-difference), respectively. $\chi^{-1}(\beta,p)$ is the chi-squared inverse function. $\mathcal{N} (0,\sigma^2)$ is the zero mean normal distribution. $\{ \emptyset \}$ denotes the empty set. The support function of $U$ at $\eta$, $h_U(\eta)$, is equal to $\sup_{u \in U} \eta^T u$.
\section{Model Descriptions}
\subsection{Plant Model $f$}
The 6dof spacecraft dynamics model is the following 
\begin{equation}
\dot{\xi}(t) = f(\xi(t),u(t))
\end{equation}
with
\begin{gather}
\begin{bmatrix}
\dot{\mathbf{r}}_a\\
\dot{\mathbf{v}}_a\\
\dot{\boldsymbol{\Theta}}\\
\dot{\boldsymbol{\omega}}_b\\
\end{bmatrix} 
= 
\begin{bmatrix}
\mathbf{v}_a\\
\mathbf{u}+\mathbf{F}-2\mathbf{\Omega} \times \mathbf{v}_a - \mathbf{\Omega} \times (\mathbf{\Omega} \times \mathbf{r}_a)\ \\
\mathbf{B}^{-1} \boldsymbol{\omega}_{b} \\
\mathbf{J}^{-1}(\mathbf{M}-\boldsymbol{\omega}_b \times \mathbf{J}\boldsymbol{\omega}_b)
\end{bmatrix}
\end{gather}
where
\begin{equation}
\xi=[\mathbf{r}_a^T\ {\mathbf{v}_a}^T\ \boldsymbol{\Theta}^T\ \boldsymbol{\omega}_b^T]^T
\end{equation}
\begin{equation}
u=[\mathbf{u}^T\ \mathbf{M}^T]^T
\end{equation}
The translational component of the model is taken from [8], and rotational component based on an example in [9]. $\xi$ is the state vector. The spacecraft is modeled as a point mass with the only external force $\mathbf{F}$. The control $u$ is such that each degree of freedom has a control channel. The translational and attitude dynamics in (2) are decoupled. $\mathbf{B} = [R_{21}R_{1i}\mathbf{e}_i \ R_{1i}\mathbf{e}_j\ \ \mathbf{e}_k]$ where
\begin{gather*}
[\mathbf{e}_{i}\ \mathbf{e}_{j}\ \mathbf{e}_{k}] =
\begin{bmatrix}
1 & 0 & 0\\
0 & 1 & 0\\
0 & 0 & 1\\
\end{bmatrix} \\
R_{ia} = \begin{bmatrix}
\cos(n_i t) & -\sin(n_i t) & 0\\
\sin(n_i t) & \cos(n_i t) & 0\\
0 & 0 & 1\\
\end{bmatrix} 
R_{1i} = \begin{bmatrix}
\cos(\psi) & -\sin(\psi) & 0\\
\sin(\psi) & \cos(\psi) & 0\\
0 & 0 & 1\\
\end{bmatrix} \\
R_{21} = \begin{bmatrix}
\cos(\theta) & 0 & \sin(\theta)\\
0 & 1 & 0\\
-\sin(\theta) & 0 & \cos(\theta)\\
\end{bmatrix} 
R_{b2} = \begin{bmatrix}
1 & 0 & 0\\
0\ & \cos(\phi) & -\sin(\phi) \\
0\ & \sin(\phi) & \cos(\phi) \\
\end{bmatrix}
\end{gather*}
\subsection{Measurement Model $h$}
The different types of measurements employed can be summarized in three main classes: 1: The straight line distance from the spacecraft to asteroid surface marker, known hereafter as a feature point, taken from lidar. 2: The location of feature points in the spacecraft body z-axis aligned camera field of view. 3: Direct angular velocity measurement $\boldsymbol{\omega}_b$, taken from IMU. The measurement model is designed to become inaccurate as the spacecraft body z-axis becomes misaligned with the feature points. \\
\begin{equation}
n(t) = [{\textbf{\emph{N}}_1(t)}^T\ ...\ {\textbf{\emph{N}}_4}(t)^T\ \emph{N}_c(t)\ {\textbf{\emph{N}}_b}(t)^T]^T
\end{equation}
is the nosie vector with $\textbf{\emph{N}}_1,...,\textbf{\emph{N}}_4,\textbf{\emph{N}}_b \in \mathbb{R}^3$. $n(t)$ is a Gaussian random vector with zero mean and constant covariance $P$. Define 
$\mathbf{d}_a^k = \mathbf{r}_a-\mathbf{p}_a^k$
where $\mathbf{p}_a^k$ is the $k$th feature point location in the asteroid fixed frame and is known. Then $\mathbf{d}_b^k = R_{ba}\mathbf{d}_a^k$. The estimated quantity $\hat{\mathbf{d}}_b^k$ is
\begin{align*}
\hat{\mathbf{d}}_b^k = \mathbf{d}_b^k + \bigg(1-\frac{-\mathbf{d}_b^k}{||\mathbf{d}_b^k||}\cdot \mathbf{e}_{k}\bigg)\textbf{\emph{N}}_k
\end{align*}
The quantity in parenthesis represents the following notion:
the further off of the center of the field of view a feature point would register on the camera, the less ability there is to accurately assess and gimbal a range sensor to the feature point. From lidar, $d_b^k = {||\hat{\mathbf{d}}_{b}^k||}_{2}$. Additionally, using the pinhole camera model provides a relation between the pixel location of the feature point on the spacecraft body z aligned camera, $\mathbf{c}=[c_1 \ c_2]^T$, to the spatial location  $\mathbf{d}_{b}$, with $f_{len}$ equal to the focal length of the camera. The measured quantity $\mathbf{c}^k$ is
\begin{align*}
\mathbf{c}^k = \bigg(\frac{f_{len}}{[0\ 0\ 1]\hat{\mathbf{d}}_{b}^k}+\emph{N}_c \bigg) \begin{bmatrix}
1 \ 0 \ 0\\
0 \ 1 \ 0\\
\end{bmatrix} \hat{\mathbf{d}}_{b}^k
\end{align*}
The last class of measurement is 
$\hat{\boldsymbol{\omega}_b} = \boldsymbol{\omega}_b + \textbf{\emph{N}}_{b}$. $f_{len}$ is assumed to hold constant for a period of time. In summary, 
\begin{gather}
y(t) = h(t,\xi(t),n(t)) \\
y = \big[d_b^1\ ...\ d_b^4\ {\mathbf{c}^1}^T\ ...\ {\mathbf{c}^4}^T\ {\hat{\boldsymbol{\omega}_b}}^T\big]^T
\end{gather}
\subsection{Navigation and Estimation}
A discrete Extended Kalman Filter is implemented with non-additive noise under the standard formulation to provide the estimation $\hat{\xi}$. The linearized state information used is given from:
\begin{gather*}
\noindent\begin{minipage}{0.16\textwidth}
\begin{equation}
A_{t} = \frac{\delta \hat{f}(t)}{\delta \xi(t)} \bigg|_{t,\hat{\xi}(t),u(t)} \label{eqab}
\end{equation}
\end{minipage}
\begin{minipage}{0.16\textwidth}
\begin{equation}
H_{t} = \frac{\delta h(t)}{\delta \xi(t)} \bigg|_{\substack{t,\hat{\xi}(t) \\
n(t)=0}} \label{eqab}
\end{equation}
\end{minipage}
\begin{minipage}{0.16\textwidth}
\begin{equation}
V_{t} = \frac{\delta h(t)}{\delta N(t)} \bigg|_{\substack{t,\hat{\xi}(t),\\
n(t)=0}} \label{eqab}
\end{equation}
\end{minipage}
\end{gather*}
The state error covariance after the filter is run at time $t$ is $\Sigma_{t}$, with $\hat{f}$ describing the spacecraft dynamics model with $\hat{\mathbf{F}}$ in place of $\mathbf{F}$.
\section{Control Problem}
$B_t$ is the linearization of $\hat{f}$ with respect to $u$ at time $t$ and $G_t$ is purposed to account for the difference between control and plant model, and is assigned before simulation. The discretized linearized stochastic open loop system is then 
\begin{gather}
 \xi (t+1) = A_t \xi (t)+B_t u (t) + G_t n (t) \\
y (t) = H_t \xi (t) +V_t n (t)
\end{gather}
\subsection{Constraints}
The constraints to keep the feature points in field of view, and box constraints on the state, respectively, are
\begin{gather}
||H_{t,fov}(\xi_{k|t}-\xi_{0|t})+h_{fov}(\xi_{0|t})||_{\infty} \leq s'_{fov} \\
||\xi_{k|t}||_{\infty} \leq s'_{\xi}
\end{gather}
which can be rewritten as $S [H_t^T \ I]^T \xi_{k|t} \leq s$, and $d$ is the number of output constraints, $d_u$ the number of control constraints, and $s=[s_{fov}^T \ s_{\xi}^T]^T$ with subscripts $fov$ denoting the rows of $H$ corresponding to outputs $\mathbf{c}^1 ... \mathbf{c}^4$, and subscript $\xi$ those applicable to box state constraints.
\subsection{Convexification}
An accurate convexification strategy is needed to perform LQMPC due to the high nonlinearity of the measurement model seen in situations where the spacecraft descends close to the feature points. Successive convexification approaches such as those in [11] can admit solutions to opimality of a non-convex optimal control (OCP) problem. Such a non-convex OCP can be constructed with $f$ propogated dynamics as opposed to (8-9), however these types of algorithms take many iterations and Jacobian evaluations per time instant. Instead of the approach of solving the non-convex OCP, the constraints introduced aim to confine the system (11-12) such that the LQMPC law is acting on a suitably convex region of the system. 
\subsection{Stochasticity}
A stochastic MPC law with horizon $N$ has tightened constraints compared to the disturbance free case as a means of probabalistically constraining the stochastic system. For a set-valued constraint $Y$, this can be written with the P-difference $Y \ominus \mathcal{Y}^{\beta}$, with $\mathcal{Y}^{\beta}$ a confidence ellipsoid with some confidence parameter $\beta$. Control strategies applied to stochastic systems to treat constraints like reference governors [7] often use a prestabilized plant via a state feedback law to propogate estimate error uncertainty that approaches a steady state value with increase in $N$. In the reference governor approach of [7] there is no consideration of control constraints, which can render a state feedback law infeasible, whereas in the case of stochastic MPC an auxiliary decision variable $v$ can be used to achieve a control objective [10]. The control $u$ is $u=Ke+v$, where $e$ is the error between (11) and its disturbance free case used in the MPC problem, $\xi_{k+1|t} = A_{t}\xi_{k|t}+B_{t}u_{k|t}$, with $\hat{\xi}(t) = \xi_{0|t}$. At time $t$ and step $k$ in the MPC horizon, $e_{k|t} = \hat{\xi}(k+t)-\xi_{k|t}$. The closed-loop component $Ke$ is used with the confidence ellipsoid to formulate the tightened constraints.

$\Phi_t$ is the stabilized error state matrix $\Phi_t = A_t+B_t K_t$ and is strictly Schur. The system (11) can then be represented with a primal dynamics part (15) and an error dynamics part (16)
\begin{align}
\zeta_{k+1|t} = A_t \zeta_{k|t} +B_t v_{k|t} \\
e_{k+1|t} = \Phi_t e_{k|t} + G_t n(k+t)
\end{align}
with $\xi_{k|t} = \zeta_{k|t}+e_{k|t}$. Consider the initial estimation error covariance $\Xi_{0|t}$ at each time taken as $\Xi_{0|t} = \Sigma_{t}$ and corresponding output error covariance $\Upsilon_{0|t}$. Over the MPC horizon they are propogated according to  
\begin{gather}
\Xi_{k+1|t} = \Phi_t \Xi_{k|t} \Phi_t^T + G_t P G_t^T \\
\Upsilon_{k|t} = H_t \Xi_{k|t} H_t^T + V_t P V_t^T
\end{gather}
with $e_{k|t} \sim \mathcal{N}(0,\Xi_{k|t})$ and similarly
$\tilde{y}_{k|t} \sim \mathcal{N}(0,\Upsilon_{k|t})$, where $\tilde{y}_{k|t} = \hat{y}(k+t)-y_{k|t}$, which is the same estimation error uncertainty propogation as [7]. The strictly Schur property ensures $\Xi_{k|t}$ does not tend to infinity with increasing $N$, and allows for a choice of $\beta$ that does not reflect less confidence as $N \rightarrow \infty$ to accomodate the increasing uncertainty. 

The output constraint at time $t$ is given as  $\zeta_{k|t} \in \tilde{Y}_{N}^{\beta} = \cap_{k=0}^{N} Y_k ^{\beta}$, where $Y_k^\beta$ is described in [7] for a prescribed $\beta$ under the assumption of a polytopic constraint set, and is given below, alongside the analogous input constraint at time $t$, 
$v_{k|t} \in \tilde{U}_{N}^{\beta} = \cap_{k=0}^{N} U_k ^{\beta}$,
\begin{multline}
Y \ominus \mathcal{Y}_k^{\beta} = Y_k^{\beta} = \\
 \{\xi \ | \ S^T_{i} [H_t^T \ I]^T \xi \leq  = s_i - h_{\mathcal{Y}_k^{\beta}}(S_i) , \ i = 1,...,d \} \\
h_{\mathcal{Y}_k^{\beta}}(S_i) = \begin{cases} \sqrt{F^{-1}(\beta,p)}\sqrt{S_i^T\ \Upsilon_{k|t} S_i} &\mbox{if } s_i \in s_{fov} \\
\sqrt{F^{-1}(\beta,p)}\sqrt{S_i^T\ \Xi_{k|t} S_i} & \mbox{if } s_i \in s_\xi \end{cases}
\end{multline} 
\begin{multline}
U \ominus \mathcal{U}_k^{\beta} = U_k^{\beta} = \\
 \{u \ | \ M^T_{i} u \leq  m_i - h_{\mathcal{U}_k^{\beta}}(M_i) , \ i = 1,...,d_u \} \\
h_{\mathcal{U}_k^{\beta}}(M_i) = \sqrt{F^{-1}(\beta,p)}\sqrt{M_{i}^T\ K_t \Xi_{k|t} K_t^T M_{i}}
\end{multline}
$Y = \{\xi \ | \ S_i^T [H_t^T \ I]^T \xi \leq  s_i, \ i = 1,...,d \}$ is the feasible set of $\xi$ given the $d$ constraints (13-14), and $S^T_i$ the $i$th row of $S$. With a slight abuse of notation $s_i \in s_{fov}$ denotes $s_i$ is found in the $s_{fov}$ partition. $\mathcal{Y}_k^{\beta} = \{ v \in \textrm{Im} \ \Upsilon_{k|t} \ | \ v^T \Upsilon_{k|t}^{+} v \leq \chi^{-1}(\beta,p) \}$ and is the confidence ellipsoid with confidence $\beta$ [7]. $h_{\mathcal{Y}_k^{\beta}}(S_i)$ is the support function of $\mathcal{Y}_k^{\beta}$ at $S_i$. $U = \{u \ | \ M_{i}^T u \leq  m_i, \ i = 1,...,d_u \}$. 

A specific terminal constraint and terminal constraint tightening to address stability of the system used in the MPC formulation is not considered in this work, although the MPC formulation that follows tends the system towards equilibrium at each time $t$.
\subsection{MPC Forumulation}
The MPC problem is reducable to a quadratic program allowing for fast simulation results. Redundant constraints resulting from (19,20) are also removed. MPC solves for the open loop control component $v$. The formulation is given as
\begin{subequations}
\begin{gather}
\min_{\delta v_t,\epsilon} J(\delta \zeta_{0|t}) = \sum_{k=0}^{N-1}(||\delta \zeta_{k|t}||_Q^2 + ||\delta v_{k|t}||_R^2) +||\delta \zeta_{N|t}||_Q^2 + {||\epsilon||}_W^2 \\
s.t. \ \delta \zeta_{k+1|t} = A_t \delta \zeta_{k|t}+B_t \delta v_{k|t} \\
S [H_{t}^T \ I]^T \zeta_{k|t} \in Y^{MPC} + \epsilon\\
M v_{k|t} \in U^{MPC}
\end{gather}
\end{subequations}
with $\epsilon$ slack variable to maintain feasibility such as in the case of empty $Y^{MPC}$, $U$ box constraints on the control channels with $ m = [m_{trans}^T \ m_{rot}^T]^T$, $Q,R,W$ weights, and 
\begin{gather*} 
Y^{MPC} =  \begin{bmatrix}
s_1 -  h_{\mathcal{Y}_N^{\beta}}(S_1) \\
\vdots \\
s_d - h_{\mathcal{Y}_N^{\beta}}(S_d)
\end{bmatrix}, \
U^{MPC} =  \begin{bmatrix}
m_1 -  h_{\mathcal{U}_N^{\beta}}(M_1) \\
\vdots \\
m_{d_u} - h_{\mathcal{U}_N^{\beta}}(M_{d_u})
\end{bmatrix}
\end{gather*}
where (21c, 21d) encode $y_{k|t} \in \tilde{Y}_{N}^{\beta}$ and $u_{k|t} \in \tilde{U}_{N}^{\beta}$ respectively. 

$\delta \zeta_{k|t} = \zeta_{k|t}-\xi^{ref}_t$, $\delta v_{k|t} = v_{k|t}-u^{ref}_t$, with $(\xi^{ref}_t,u^{ref}_t)$ the equilibrium pair for the system (21b), with $\xi_t^{ref}$ specified before the computation of (21). $u^{ref}_t$ is computed as $u^{ref}_t = (B_t^T B_t)^{-1} B_t (I-A_t)\xi^{ref}_t$.
\section{Spacecraft Descent Simulation and Results}
Two trajectories are simulated and discussed in this section. Leg 1 was chosen such that the output dynamics are highly nonlinear as the spacecraft descends and the feature points move further out of the field of view. Leg 2 simulates a landing approach ending around 10 meters above the surface. Both trajectories consist of a constant descent rate phase and hover at a point phase. At each time step the linearizations (8-10) are updated by the EKF which are also used in the MPC law. A desired position reference $\mathbf{r}_{a}^{ref}$ is first used to calculate $\boldsymbol{\Theta}^{ref}$ before the simulation starts, so as to center the camera to the average location of all feature points $\mathbf{d}_a^{avg}$ (per time step), and then the complete reference $\xi^{ref} = [\mathbf{r}_a^{{ref}^T} \mathbf{0}^T \boldsymbol{\Theta}^{{ref}^T} \mathbf{0}^T]^T$ is fed to the simulation. The requisite Euler angles are extracted from the rotation matrix $R_{ba} = R_{b2}R_{21}R_{1i}R_{ia}$. Due to the nonunique properties of Euler angles for a given rotation, the same seeding used to initialize the rotation matrix in $\boldsymbol{\Theta}_0^{ref}$ is kept throughout. $\boldsymbol{\Theta}_0^{ref}$ is initialized $R_{ba}\mathbf{r}_a (0) = ||\mathbf{d}_a^{avg}||_2 \ [0\ 0\ 1]^T$. $f_{len}$ remains constant throughout the trajectory, and is initialized as $\mathbf{d}_b^k$ that is best projected along the spacecraft body z axis. All plots and code implemented in MATLAB R2018b.

{\centering
\includegraphics[scale=0.28]{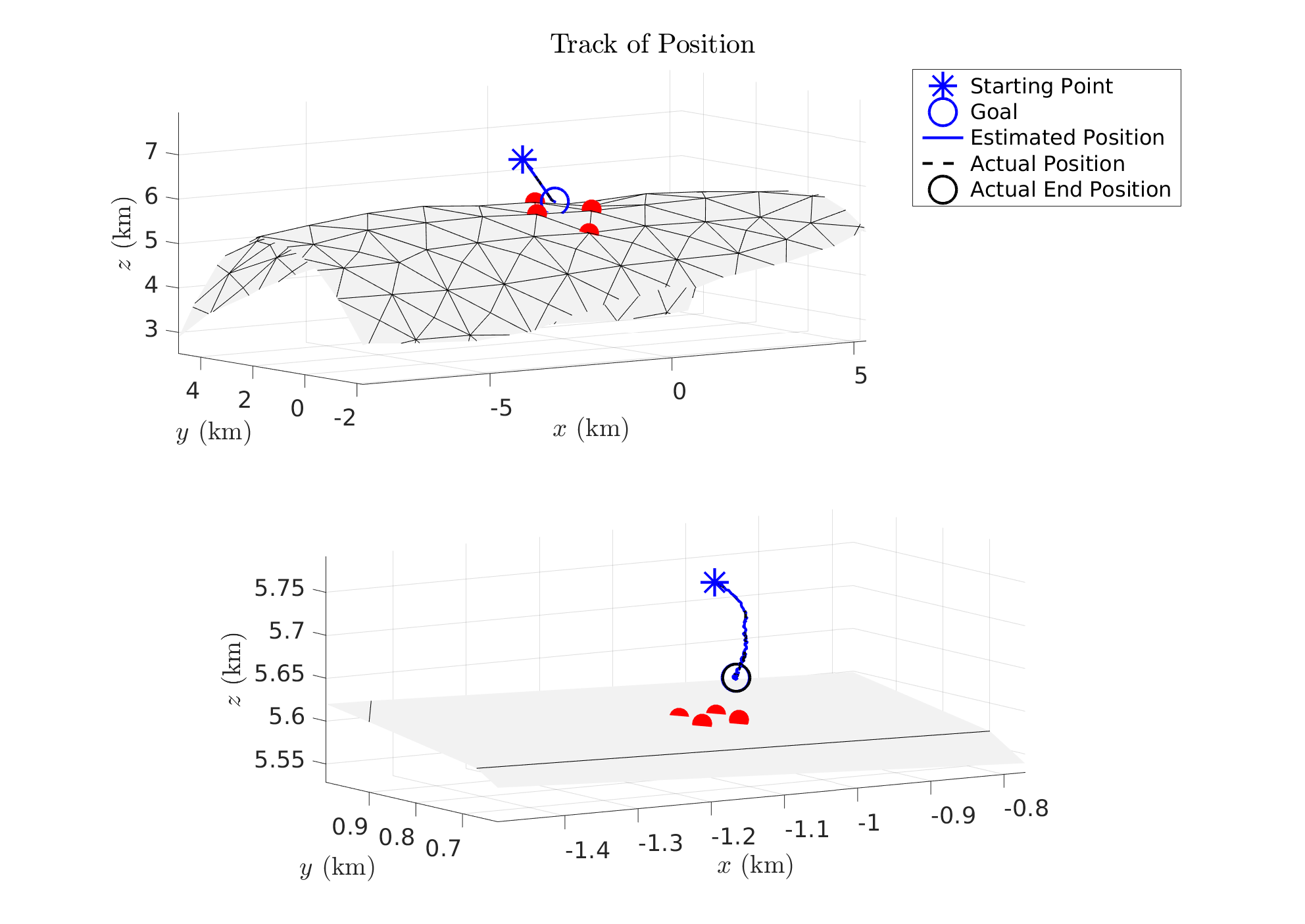} 
\captionof{figure}{\footnotesize The upper plot shows leg 1 that begins and ends further away from the asteroid surface, as compared to leg 2 in the lower plot. The corresponding feature points are also plotted in red. The asteroid is in grey. Constraint satisfaction are prioritized over tracking, as seen in the lower plot.} 
}
{\centering
\includegraphics[scale=0.28]{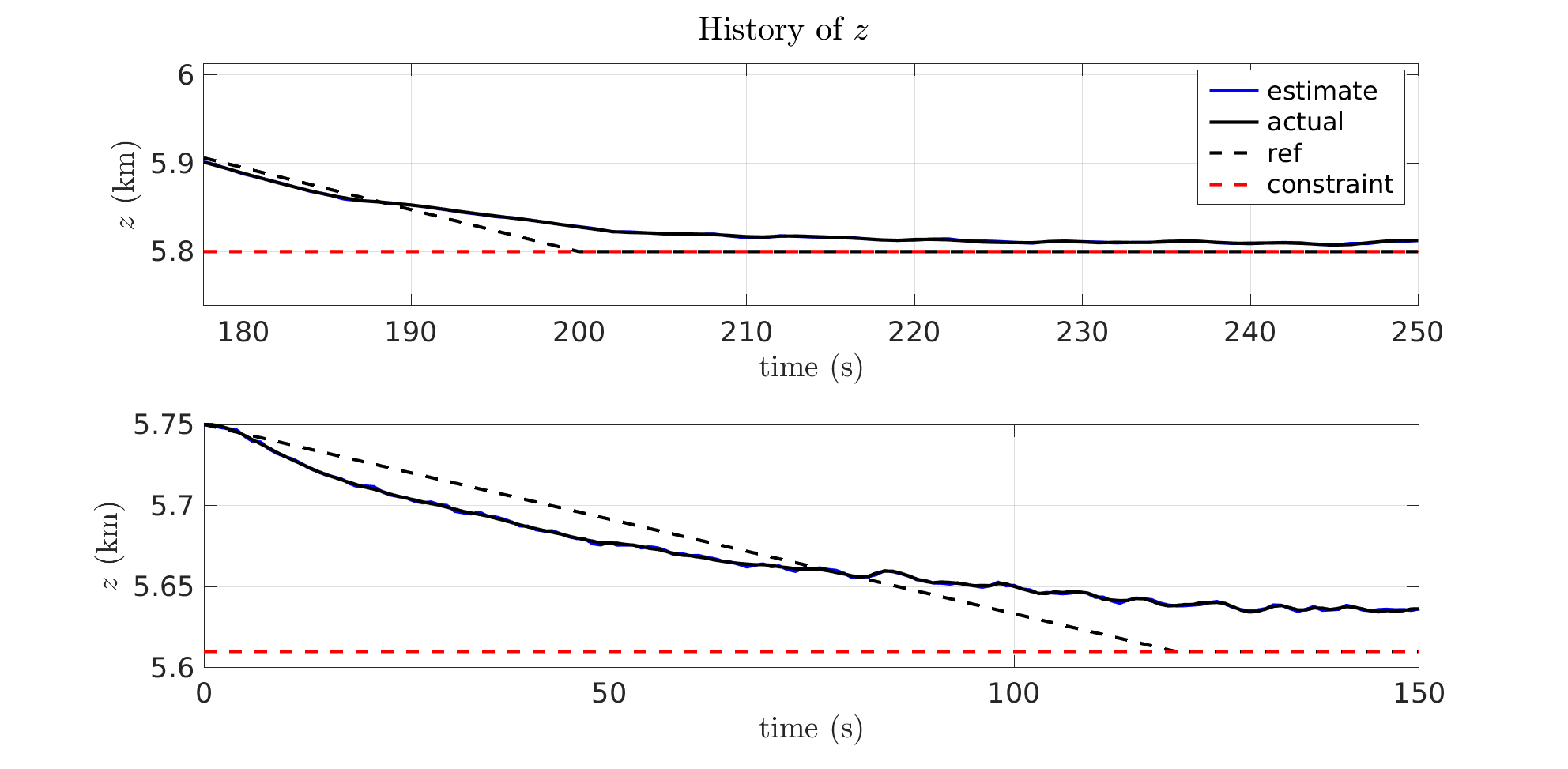}
\captionof{figure}{\footnotesize Constraints from (14) are plotted. In both legs the field of view constraint (13) is encountered before (14). The untightened constraints are shown.}
} 
{\centering
\includegraphics[scale=0.28]{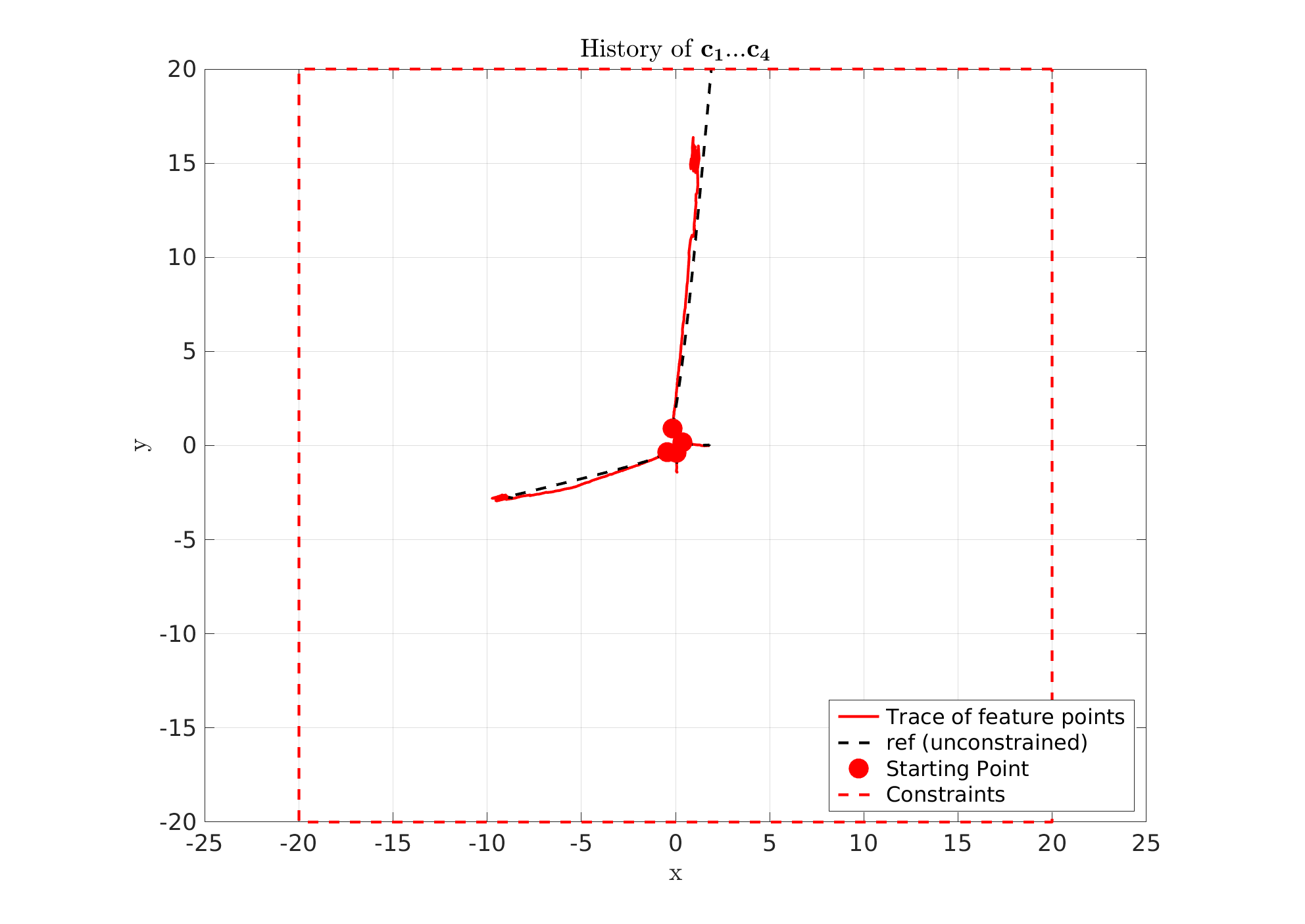} 
\captionof{figure}{\footnotesize Constraints from (13) of leg 1 are plotted. For this trajectory in the unconstrained case some feature points locations are expected to expand quickly out of view upon further descent. The untightened constraints are shown, along with the feature point trace from the reference trajectory which extends beyond the constraints.}
}
{\centering
\includegraphics[scale=0.28]{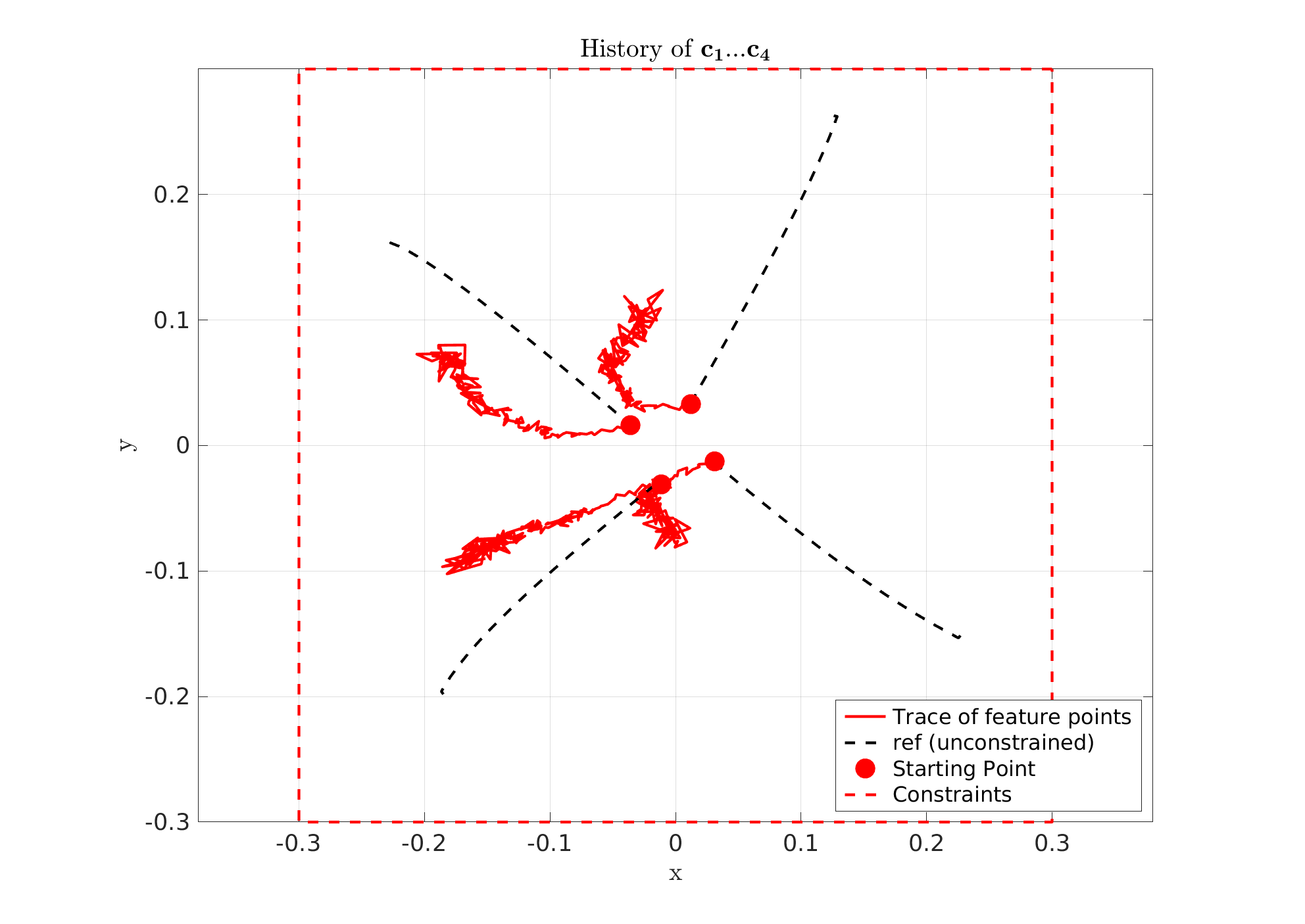} 
\captionof{figure}{\footnotesize Constraints from (13) of leg 2 are plotted. The untightened constraints are shown, along with the feature point trace from the reference trajectory.}
}
 
The simulation parameters are given below. The $\texttt{place}$ command in MATLAB was used to compute $K_t$. The placed eigenvalues are the same across all times for both legs.

\begin{center}
\scalebox{0.85}{
 \begin{tabular}{c c c} 
 \hline
 \hline
 & Leg 1 & Leg 2 \\ [0.5ex] 
 Discretization $(sec)$ & 1 & 1 \\ [0.5ex]
 $\beta$ & 0.95 & 0.95 \\ [0.5ex]
 $N$ & 20 & 20 \\ [0.5ex]
 $s'_{fov}$ & 20 & 0.3 \\ [0.5ex]
 $s'_{\xi}= z_{cnstr}$ & 5.8 & 5.61 \\ [0.5ex]
 $m_{trans} \ (km/s^2)$ & 0.002 & 0.002 \\ [0.5ex]
 $m_{rot}$ & 0.5 & 0.5 \\ [0.5ex]
 $\mathbf{r}_0^{ref}$ & [-1.57 1.32 6.75]$^T$ & [-1.07 0.82 5.75]$^T$ \\ [0.5ex]
 $\mathbf{r}_{end}^{ref}$ & [-1.07 0.82 5.8]$^T$ & [-1.07 0.82 5.61]$^T$ \\ [0.5ex]
 $Q$ & $\left[ \begin{smallmatrix} 0.1 & & & \\  & 0.0001 & & \\ & & 0.1 & \\ & & & 0.0001 \end{smallmatrix} \right] $ & $\left[ \begin{smallmatrix} 10 & & & \\  & 0.0001 & & \\ & & 10 & \\ & & & 0.0001 \end{smallmatrix} \right] $ \\ [0.5ex]
 $R$ & $\left[ \begin{smallmatrix} 50 & & & \\  & 50 & & \\ & & 1000 & \\ & & & 0.001 \end{smallmatrix} \right] $ & $\left[ \begin{smallmatrix} 5000 & & & \\  & 5000 & & \\ & & 5000 & \\ & & & 1 \end{smallmatrix} \right] $ \\ [0.5ex]
 $W$ & $\left[ \begin{smallmatrix} 10^5 \end{smallmatrix} \right] $ & $\left[ \begin{smallmatrix} 10^9 \end{smallmatrix} \right] $ \\ [0.5ex]
 $P$ & \multicolumn{2}{c}{diag([0.001]$^2$)}\\ [0.5ex]
 $G$ & \multicolumn{2}{c}{0.00001$J_1$}\\
 \hline
 \hline
\end{tabular}}
\captionof{table}{\footnotesize Table of selected parameters. The only state constraint is on $z$ and is $z_{cnstr}$. The elements in the constraint inequalities $s$ and $m$ are equal in absolute value for those in the same subscript, thus a single element is reported here. $Q,R,W$ are block diagonal, again with a single element in each block reported, with $Q$ partitoned along $\mathbf{r}_a,\mathbf{v}_a,\boldsymbol{\Theta},\boldsymbol{\omega}_b$, and $R$ along $u_x,u_y,u_z,\mathbf{M}$. $J_1$ is the matrix of ones.}
\end{center}

\section{Conclusion}
In this paper stochastic tube LQMPC was implemented on a spacecraft mission to descend to an orbiting asteroid surface in the presense of uncertainty in gravity and in the measurement model. A technology based measurement model was first developed; this led to a camera field of view constraint. Two trajectories were simulated with focus on the trace of feature point locations in the camera. Future works of interest concerning the control are formulating an OCP with terminal constraints constructed to address recusive feasability and stability, as well as addressing the degree of convex approximation used in the MPC law and using a non-convex OCP. A coupled translational and rotational actuation model can also be handled by the processes laid out in this paper. A final direction of future work is a scheme to estimate $G$, the gravity model mismatch.

\end{document}